\RequirePackage{fix-cm}
\documentclass[smallcondensed,envcountsect,envcountsame]{svjour3} 
\smartqed  
\usepackage{graphicx}
\usepackage{amsmath,amssymb}

\usepackage[all]{xy}
\newenvironment{spmatrix}{\left ( \begin{smallmatrix}} {\end{smallmatrix}\right )}
\newenvironment{sbmatrix}{\left [ \begin{smallmatrix}} {\end{smallmatrix}\right ]}

\spnewtheorem*{pf}{Proof.}{\it}{\rm}
\spnewtheorem{rem}[theorem]{Remark}{\bf}{\rm}
\spnewtheorem{exm}[theorem]{Example}{\bf}{\rm}
\renewcommand*{\proofname}{Proof.}
{\rm}
\spnewtheorem*{THM}{Theorem}{\bf}{\it}

\begin{document}

\title{A characterization of Veech groups in terms of origamis
}


\author{Shun Kumagai
}


\institute{S. Kumagai \at
              Research Center for Pure and Applied Mathematics, 
              Graduate School of Information Sciences, 
              Tohoku University, Sendai 980-8579, Japan.  \\
              \email{shun.kumagai.p5@dc.tohoku.ac.jp}    \\
              This work was supported by Japan Student Service Organization. 
}

\date{Received: date / Accepted: date}

\maketitle

\begin{abstract}
Schmith\"usen proved in 2004 that the Veech group of an origami is closely related to a subgroup of the automorphism group of the free group $F_2$. 
This result is significant in the sense that the framework of approachable Veech groups is greatly extended. 
In this paper, we continue the analysis and consider what kind of settings of flat surfaces allow Veech groups to be characterized combinatorially like origamis.
We show that elements in the Veech group of a flat surface with two finite Jenkins-Strebel directions are characterized to allow a concurrence between two `origamis' defined by geodesics in the surface. 
In the proof we use an observation presented by Earle and Gardiner that a flat surface with two finite Jenkins-Strebel directions is decomposed into a finite number of parallelograms and is proved to be of finite analytic type. 
Using our results we can decide whether a matrix belongs to the Veech group for various kinds of flat surfaces of finite analytic type. 

\end{abstract}

\section{Introduction}
\label{intro}
In this paper we consider Riemann surfaces of finite analytic type. 
A pair $(R,\phi)$ of a Riemann surface $R$ and a holomorphic quadratic differential $\phi$ on $R$ is called a flat surface \cite{V} (also called a half-translation surface). 
For a flat surface we can define the notion of local `affine' geometry including a flat metric, straight lines, directions, etc. 
This notion appears in the Teichm\"uller theorem which says that every point in the Teichm\"uller space $\mathcal{T}(R)$ is represented by an affine deformation with minimal dilatation of the base point equipped with a holomorphic quadratic differential $\phi$ on $R$.  
The collection of all affine deformations of a flat surface $(R,\phi)$ gives an isometric embedding $\mathbb{H}\hookrightarrow \mathcal{T}(R)$ of the upper half plane $\mathbb{H}$ and its image $\Delta(R,\phi)\subset \mathcal{T}(R)$, called the Teichm\"uller disk.

 
 The Veech group $\Gamma(R,\phi)$ of a flat surface $(R,\phi)$ is the group of derivatives of all self affine deformations of $(R,\phi)$. 
This group is a discrete subgroup of the group $PSL(2,\mathbb{R})$ acting on the upper half plane $\mathbb{H}$ as a group of M\"obius transformations. There is a correspondence between this model and the projected image $C(R,\phi)$ of $\Delta(R,\phi)$ in the moduli space $M(R)$ of Riemann surfaces via the natural projection $\mathcal{T}(R)\rightarrow M(R)$. 
Veech groups originally have been studied in the context of billiard dynamics and the geodesic flow on billiard tables
 \cite{V}. 

As necessary we consider a `square-tiled' flat surface, a concept that comes from finite copies of the Euclidian unit square and their natural coordinates. Such a surface is identified by a combinatorial data representing how squares adjoin each other, we call this data an origami. 
In this case each point in $\Delta(R,\phi)$ with a parameter $\tau\in\mathbb{H}$ can be seen as the surface given by deforming the Euclidian unit square to a parallelogram of modulus $\tau$. 

Schmith\"usen \cite{S1} proved that the Veech group of an origami is closely related to a subgroup of the automorphism group of the free group $F_2$. This observation gives a computable characterization of the Veech group as a subgroup of  the group $PSL(2,\mathbb{Z})$ of finite index. 
So for any origami, $C(R,\phi)$ is an algebraic curve defined over the algebraic number field $\bar{\mathbb{Q}}$, called an origami curve. 
The action of the absolute Galois group $Gal(\bar{\mathbb{Q}}/\mathbb{Q})$ on origami curves are observed by M\"oller \cite{M1}.
He proved that origami curves form geometric components of a Hurwitz space defined over $\bar{\mathbb{Q}}$ and that the Galois actions on origami curves and origamis are `compatible'. 
With this result he also gave an another proof of the inclusion of $Gal(\bar{\mathbb{Q}}/\mathbb{Q})$ in the Grothendieck-Teichm\"uller group. Origamis have been studied in the context of Teichm\"uller theory and number theory. (See for instance \cite{H}, \cite{HS1}, \cite{HS2},  \cite{L}, and \cite{M2}.)

In this paper we approach flat surfaces from general settings. We consider what kind of settings allow Veech groups to be characterized combinatorially like origamis. 
In section \ref{sec3-1} we see that two finite Jenkins-Strebel directions $\theta_1,\theta_2$ of a flat surface $(R,\phi)$ induce a unique decomposition of $R$ into finite numbers of parallelograms. This observation was already discussed in \cite{EG}. 
As we will state in section \ref{F2} there is an action of $F_2$ to signed parallelograms defined by geodesics in the directions $\theta_1,\theta_2$. This action is identified with an origami. 

In section \ref{sec3-3} we will see that a flat surface $(R,\phi)$ with two finite Jenkins-Strebel directions $(\theta_1,\theta_2)$ is uniquely determined by the decomposition into parallelograms, called the \textit{P-decomposition} $P(R,\phi,(\theta_1,\theta_2))$, which consists of data of angles, a list of moduli, and an origami. 
Variations of a P-decomposition under affine deformations are observed as the variations in the plane
, which can be seen as a free action of derivatives of self affine deformations. 
As a result the comparison between the initial P-decomposition and the terminal decomposition determines the existence of an affine map with given derivative. 



For a flat surface $(R,\phi)$ we denote the set of Jenkins-Strebel directions by $J(R,\phi)$. 
We define an action $A\in PSL(2,\mathbb{R})$ on angles, parallelograms and P-decompositions by the variation of them in the plane under an affine map with derivative $\pm A$. 
With these notations the main result is stated as follows. 
\begin{THM}
Let $(R,\phi)$ be a flat surface of finite analytic type with two distinct Jenkins-Strebel directions $\theta_1,\theta_2\in J(R,\phi)$. ${A}\in PSL(2,\mathbb{R})$ belongs to $\Gamma (R,\phi)$ if and only if $A\theta_1,A\theta_2$ belongs to $ J(R,\phi)$ and $A\cdot P(R,\phi,(\theta_1,\theta_2))$ is isomorphic to $P(R,\phi,(A\theta_1,A\theta_2))$. 

\end{THM}

\section{Preliminaries}
\label{sec:2}

\subsection{Flat structures and Veech groups}
\label{sec:2-1}
Let $R$ be a Riemann surface of finite analytic type $(g,n)$ with $3g-3+n>0$. 
\\
First we give some definitions related with the \textit{flat structures} \cite{V} on Riemann surfaces, which are always defined by their \textit{quadratic differentials}. We present some definitions and properties describing `affine deformations of a flat structure'. For details, see \cite{EG}, \cite{GJ} for instance. 


\begin{definition}Let $R_i$ $(i=1,2)$ be Riemann surfaces homeomorphic to $R$. 
\begin{enumerate}
\item A homeomorphism $f:R\rightarrow R_1$ is \textit{quasiconformal mapping} if $f$ has locally integrable distributional derivatives $f_z,f_{\bar{z}}$ and there exists $k<1$ such that $|f_{\bar{z}}|\leq k|f_z|$ holds almost everywhere. 
We denote the group of quasiconformal mappings of $R$ onto itself by $QC(R)$. 

\item We say two quasiconformal mapping $f_i:R\rightarrow R_i$ $(i=1,2)$ are \textit{Teichm\"uller  equivalent} if there is a conformal map $h:R_1\rightarrow R_2$ homotopic to $f_2\circ f_1^{-1}:R_1\rightarrow R_2$. 

\item We define the \textit{Teichm\"uller space} $\mathcal{T}(R)$ of Riemann surface $R$ as the space of Teichm\"uller equivalence classes of quasiconformal mappings from $R$. We define the \textit{mapping class group} $\mathrm{Mod}(R)$ by the group of homotopy classes of every element in $QC(R)$.

\item For $f\in QC(R)$, we define $\rho_f:\mathcal{T}(R)\rightarrow\mathcal{T}(R)$ by $[g]\mapsto[g\circ f^{-1}]$ for every quasiconformal mapping $g$ from $R$. Now the homomorphism $QC(R)\ni f\mapsto \rho_f$ factors through $\mathrm{Mod}(R)$. We define the \textit{moduli space} $M(R)$ by the quotient $\mathcal{T}(R)/\mathrm{Mod}(R)$.



\end{enumerate}
\end{definition}

\begin{definition}
A \textit{Beltrami differential} $\mu$ on $R$ is a tensor on $R$ whose restriction to each chart $(U,z)$ on $R$ is of the form $\mu(z){d\bar{z}}/{dz}$ where $\mu$ is a measurable function on $z(U)$. 
\end{definition}
We define a norm of a Beltrami differential $\mu$ by the essential supremum of $|\mu(z)|$ where $(U,z)$ runs all charts on $R$ and $\mu$ is locally represented as above definition. A Beltrami differential whose norm is less than $1$ is called a \textit{Beltrami coefficient} and we denote the space of Beltrami coefficients on $R$ by $B(R)$. 
From the measurable Riemann mapping theorem we see that for any $\mu\in B(R)$ there exists a unique quasiconformal mapping from $R$ whose Beltrami coefficient $\mu_f$ equals to $\mu$. 
So we define $\Phi:B(R)\rightarrow \mathcal{T}(R)$ by $\Phi(\mu):=[f]$ where $\mu=\mu_f$.

\begin{definition}
A \textit{holomorphic quadratic differential} $\phi$ on $R$ is a tensor on $R$ whose restriction to each chart $(U,z)$ on $R$ is of the form $\phi(z)dz^2$ where $\phi$ is a holomorphic function on $U$. \end{definition}
Let $p_0\in R$ be a regular point of $\phi$ and $(U,z)$ be a chart around $p_0$. Then $\phi$ defines a natural coordinate ($\phi$\textit{-coordinate}) $\zeta (p)=\int^p_{p_0} \sqrt{\phi(z)}dz$ on $U$, on which $\phi=d\zeta^2$. $\phi$-coordinates give an atlas on $R^*=R\setminus \mathrm{Crit}(\phi)$ whose any coordinate transformation is of the form $\zeta \mapsto \pm \zeta +c\ (c\in \mathbb{C})$. 
Such a structure, which is a maximal atlas whose any coordinate transformation consists of half-turn and translation is called a \textit{flat structure} on $R^*$. A flat structure (or such an atlas) on $R^*$ determines a holomorphic quadratic differential $\phi=d\zeta ^2$ on $R$. We define a norm of a quadratic differential by the surface integral on $R$ with the natural coordinates and denote by $\|\cdot\|$. The space $A(R)$ of quadratic differential of finite norm on $R$ is known to be a vector space of complex dimension $3g-3+n$. In the following we take $\phi\in A(R)$ and assume $R=R^*$ by puncturing at the points of the discrete set $\mathrm{Crit}(\phi)\subset R$ if necessary. 

\begin{rem}
For a Riemann surface $R_0$ of type $(g,n)$ 
the Teichm\"uller space 
$\mathcal{T}(R_0)$ is  a complex manifold of dimension $3g-3+n$ which is homeomorphic to the unit ball in $\mathbb{C}^{3g-3+n}$. For $[R,f]\in\mathcal{T}(R_0)$, $A(R)$ is known to be the cotangent space $T_{[R,f]}\mathcal{T}(R_0)$ which is a complex vector space of dimension $3g-3+n$. 
\end{rem}

Fix a holomorphic quadratic differential $\phi$ on $R$ satisfying $\|\phi\|=1$. It gives the Beltrami differential $\bar{\phi}/|\phi|$ whose restriction to each chart $(U,z)$ is $\bar{\phi(z)}/|\phi(z)|\cdot d\bar{z}/dz$. 
For each $t\in \mathbb{D}$, $t\bar{\phi}/|\phi|\in B(R)$ and thus the map $\mathbb{D}\ni t\mapsto [f_t]\in \mathcal{T}(R)$ is well-defined. In fact this defines a holomorphic, isometric embedding with respect to the Poincar\'e metric and the Teichm{\"u}ller metric (see \cite{Ga}).  
We call this embedding $\iota_\phi : \mathbb{D}\hookrightarrow \mathcal{T}(R)$ the \textit{Teichm\"uller embedding} and its image $\Delta_\phi = \iota_\phi (\mathbb{D})$ the \textit{Teichm\"uller disk}. 


\begin{definition}Let $R,\phi$ be as above. 
\begin{enumerate}
\item For $A=\begin{spmatrix}a&b\\ c&d\end{spmatrix}\in SL(2,\mathbb{R})$, we define $T_A:\mathbb{C}\rightarrow\mathbb{C}$ by $\zeta=\xi+i\eta\mapsto (a\xi+c\eta)+i(b\xi+d\eta)$. (The derivative of $T_A$ equals $A$.)
\item A quasiconformal homeomorphism $f:R\rightarrow R$ is called an \textit{affine map} on $(R,\phi)$ if $f$ is locally affine (i.e. of the form $f(\zeta)=T_A(\zeta) +k$) with respect to the 
$\phi$-coordinates. We denote the group of all affine map on $(R,\phi)$ by $\mathrm{Aff}^+(R,\phi)$. 
\item For each $f\in \mathrm{Aff}^+(R,\phi)$, the local derivative $A$ is globally defined up to a factor $\{\pm I\}$ independent of coordinates of $u_\phi$. We call the map $D:\mathrm{Aff}^+(R,\phi)\rightarrow PSL(2,\mathbb{R}) :f\mapsto \bar{A}$ the \textit{derivative} map and its image  $\Gamma (R,\phi):=D(\mathrm{Aff}^+(X,\phi))$ the 
\textit{Veech group}.
\item Let $R$ be a Riemann surface of finite analytic type and $\phi$ be a non-zero, integrable, holomorphic quadratic differential on $R$. We call such a pair $(R,\phi)$ a \textit{flat surface}. We say that flat surfaces $(R,\phi),(S,\psi)$ are \textit{isomorphic} if there is a conformal map $f:R\rightarrow S$ such that $\phi =f^*\psi$. 
\end{enumerate}
\end{definition}

\begin{rem}


Let two flat surfaces $(R,\phi),(S,\psi)$ be isomorphic with $f:R\rightarrow S$ as above definition. Then $f$ is locally affine with respect to $\phi$-coordinates and $\psi$-coordinates with derivative $[I]\in PSL(2,\mathbb{R})$. 
\end{rem}



\begin{lemma}[{\cite[Theorem1]{EG}}]\label{mob}Let $(R,\phi)$ be a flat surface and $f\in QC(R)$. 
Then $\rho_f$ maps $\Delta_\phi$ onto itself if and only if $f$ is homotopic to an element in $\mathrm{Aff}^+(R,\phi)$. 
Furthermore in this case, $\rho_f(\Phi(t\bar{\phi}/|\phi|))=\Phi(D(f)^*(t)\bar{\phi}/|\phi|)$ for each $t\in \mathbb{D}$ where $\begin{sbmatrix}a&b\\ c&d\end{sbmatrix}^*(\tau):=\frac{-a\tau+b}{c\tau-d}$ for $\tau=${\scriptsize$\sqrt{-1}\cdot$}$\frac{1+t}{1-t}\in\mathbb{H}$. 

\end{lemma}

Recall that the Veech group $\Gamma(R,\phi)$ is the group of derivatives of elements in $\mathrm{Aff}^+(R,\phi)$. 
Lemma \ref{mob} says that how the Teichm\"uller disk $\Delta_\phi$ projects into $M(R)$ can be seen from $\Gamma(R,\phi)$ acting on $\mathbb{H}$. 
It is first observed by Veech \cite{V} that Veech group is a discrete group. The Veech group determines the projected image $C_\phi$ of $\Delta_\phi$ in $M(R)$ in the sense that $C_\phi$ is isomorphic to the mirror image of $\mathbb{H}/\Gamma (R,\phi)$ as an orbifold. If $\Gamma (R,\phi)$ has a finite covolume then $C_\phi$ can be seen as a Riemann surface of finite analytic type, called the \textit{Teichm\"uller curve} induced by $\phi$. 

\begin{rem}
For a flat surface $(R,\phi)$ if $\sqrt{\phi}$ (whose restriction to a chart $(U,z)$ is $\sqrt{\phi(z)}dz$) gives an Abelian differential on $R$ then we can take the subatlas of $u_\phi$ so that all coordinate transformation is of the form $\zeta\mapsto \zeta +c$ $(c\in\mathbb{C})$. We call such a structure \textit{translation structure} and such $\phi$ \textit{Abelian}.  

In Abelian case the derivative map is well-defined onto $SL(2,\mathbb{R})$ and the Veech group $\Gamma(R,\phi)$ is defined to be a subgroup of $SL(2,\mathbb{R})$. We denote the projected class of each $A\in SL(2,\mathbb{R})$ by $\bar{A}$ and the projected Veech group by $\bar{\Gamma}(R,\phi)$ in Abelian case. 

\end{rem}

\subsection{Origamis}
\label{sec2-2}

Origamis \cite{L} are 
combinatorial objects which induce "square-tiled" flat structures, whose Veech groups can be characterized as a projected image of a subgroup of $\mathrm{Aut}(F_2)$. 
They are good examples in the sense that they always produce Teichm\"uller curves defined over $\bar{\mathbb{Q}}$. 

They are also studied in the context of the Galois action on combinatorial objects as well as \textit{dessins d'enfants}, a crucial result is given by M\"oller \cite{M1} and some of study is described in \cite{HS1}. 

\begin{definition}
An \textit{origami} is a topological covering $p:R\rightarrow E$ from a connected oriented surface $R$ to the torus $E$ ramified at most over one point $\infty\in E$. We say two origamis $\mathcal{O}_j=(p_j:R_j\rightarrow E)$ $(j=1,2)$ are \textit{equivalent} if there exists a homeomorphism $\varphi :R_1\rightarrow R_2$ such that $p_1= p_2\circ \varphi$. 
\\We consider the branch points of $p$ as marked points. 
\end{definition}

\begin{exm}The origami shown in following figure is called \textit{L-shaped origami} $L(2,3)$. The surface $R^*$ is of type $(2,2)$ and Schmithu\"sen \cite{S1} showed that the Veech group of $L(2,3)$ is a non congruence group. 

\begin{figure}[htbp]
\begin{center}
  \includegraphics[width=90mm]{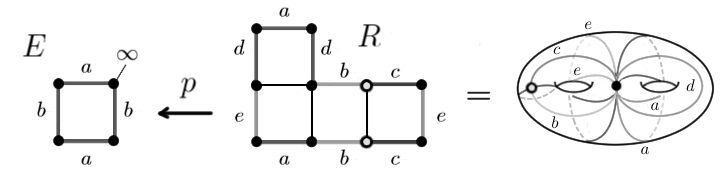}
  \caption{L-shaped origami $L(2,3)$: edges with the same letters are glued. }
\label{Fig.1}    
\end{center}
\end{figure}
In general an origami can be seen as a surface obtained by gluing finite unit squares at edges in the way that natural coordinates $z$ of squares give a globally defined Abelian differential $dz$. Such a way of gluing is called an \textit{origami rule}. 
\end{exm}
By a general theory of covering maps, we have following characterizations for an origami similar to a \text{dessin d'enfants}. See \cite{HS1} for details.

\begin{proposition}\label{origami}
An origami of degree $d$ is up to equivalence uniquely determined by each of the following. 

\begin{enumerate}
\item{An origami-rule for $d$ unit square cells. }
\item {A finite oriented graph $(\mathcal{V, E})$ such that $|\mathcal{V}|=d$ and every vertex has precisely two incoming edges and two outgoing edges, with both of them consist of edges labeled with $x$ and $y$. }
\item{A monodromy map $m:F_2 \rightarrow S_d$ up to conjugation in $S_d$. }
\item{A subgroup $H$ of $F_2$ of index $d$ up to conjugation in $F_2$. }
\end{enumerate}
\end{proposition}
Note that $F_2=\langle x,y \rangle$ is the free group generated by two elements. For each origami we can consider an action $F_2 \hspace{0pt}^\mathbb{y} \mathcal{V}$ which comes from the monodromy. 
\\
A combinatorial characterization of the Veech group for an origami is given as in the following lemma. 
\begin{lemma}[Schmith\"usen {\cite[Lemma 2.8]{S1}}]\label{diagram} Let $\mathcal{O}=(p:R\rightarrow E)$ be  an origami and $\mathbb{H}$ be the universal covering space of $E^*$ equipped with induced Abelian differential $\psi$. 
Then there exists a subgroup $\mathrm{Aut}^+(F_2)$ of $\mathrm{Aut}(F_2)$ and following exact (in horizontal direction) and commutative diagram.
\[
\xymatrix{
1\ar[r]&\mathrm{Gal}(\mathbb{H}/E^*)\ar@{^{(}->}[r]\ar[d]^{\cong}&\mathrm{Aff}^+(\mathbb{H},\psi)\ar[r]^D \ar[d]^{\cong} &\mathrm{SL}_2(\mathbb{Z})\ar[r]\ar[d]^{\cong}&1
\\
1\ar[r]&\mathrm{Inn}(F_2)\ar@{}[ru]|(0.5){\rm \circlearrowleft} \ar@{^{(}->}[r]&\mathrm{Aut}^+(F_2)\ar[r]\ar@{}[ru]|(0.5){\rm \circlearrowleft}&\mathrm{Out}^+(F_2)\ar[r]&1
}\]
Furthermore, the subgroup $\mathrm{Aff}^+(H)<\mathrm{Aut}^+(F_2)$ which corresponds to \\$\mathrm{Aff}^+(R^*,\phi)<\mathrm{Aff}^+(\mathbb{H},\psi)$ in this diagram coincides with   
$\mathrm{Stab}_{\mathrm{Aut}^+(F_2)}(H)$. 

\end{lemma}

 \begin{corollary}[{\cite[Corollary 2.9]{S1}}]\label{SL}
 $\Gamma (R^*,\phi)$ is a subgroup of $SL(2,\mathbb{Z})$ of finite index. 
 \begin{pf}
An automorphism of $F_2$ does not change an index of subgroup. So an $\mathrm{Aut}^+(F_2)$-orbit of $H$ consists of $F_2$-subgroups of common finite index. By proposition \ref{origami} we see that such subgroups are origamis of common degree and hence they are finite. 
\qed\end{pf}
 \end{corollary}

\begin{rem}
A subgroup $\Gamma<SL(2,\mathbb{Z})$ of finite index gives a finite covering $\mathbb{H}/\bar{\Gamma}\rightarrow \mathbb{H}/PSL(2,\mathbb{Z})$ unramifed over the singularities of $\mathbb{H}/PSL(2,\mathbb{Z})$. This is uniquely extended to a meromorphic function ramified over at most $i, e^{\pi i/3}, \infty \in \hat{\mathbb{C}}$.  Hence Corollary \ref{SL} implies that the translation structure for an origami always gives a Teichm\"uller curve which is a Belyi curve. The complex structure of Teichm\"uller curve of an origami can be written by its dessin d'enfants. 
\end{rem}


\section{Characterization of Veech group}
\subsection{Decompositions with the $\phi$-metric}\label{sec3-1}
Let $R$ be a Riemann surface of finite analytic type and $\phi\in A(R)$. 
\\
The Euclidian metric lifts via $\phi$-coordinates to a flat metric $g(\phi)$ on $R$. We call this metric the $\phi$\textit{-metric} and geodesics of $g(\phi)$ the $\phi$\textit{-geodesics}. The $\phi${-geodesics} are mapped to the geodesics of the complex plane (i.e. line segments) by the $\phi$-coordinates. 
Since coordinate transformations of a flat structure do not change the slope of line segment, the slopes of the $\phi$-geodesics are well-defined. For any $\alpha \in S^1$ the metrics $g(\phi), g(\alpha\phi)$ coincide and $\alpha\phi$-coordinates are products of $\phi$-coordinates and $\sqrt{\alpha}$. So any $\phi$-geodesic is a horizontal (slope $0$) $\alpha\phi$-geodesic for some $\alpha\in S^1$. 

\begin{definition}\ 

\begin{enumerate}
\item The \textit{direction} of a $\phi$-geodesic $\gamma$ is $\theta\in [0, \pi)$ where $\gamma$ is horizontal $e^{2\sqrt{-1}\theta}\phi$-geodesic. 
\item The $\phi$\textit{-cylinder} generated by a $\phi$-geodesic $\gamma$ is the union of all $\phi$-geodesics parallel (with same direction) and free homotopic to $\gamma$. We define the direction of a $\phi$-geodesic by the one of its generator. 
\item $\theta\in [0, \pi)$ is \textit{Jenkins-Strebel direction} of $(R,\phi)$ if almost every point in $R$ lies on some closed $\phi$-geodesic of direction $\theta$. 
We denote the set of Jenkins-Strebel directions by $J(R,\phi)$. 

\end{enumerate}
\end{definition}
Note that any Jenkins-Strebel direction of flat surface of finite analytic type is \textit{finite}, namely there are at most finitely many $\phi$-cylinders of that direction in $R$. For the existence of a holomorphic differential with one Jenkins-Strebel direction, the following result is known. 
\begin{proposition}[Strebel \cite{St}]
Let $\gamma=(\gamma_1,...,\gamma_p)$ be a finite `admissible' curve system on $R$, which satisfies bounded moduli condition for $\gamma$. Then for any $b=(b_1,...,b_p)\in\mathbb{R}_+^p$ {there exists $\phi\in A(R)$ such that $0$ is a Jenkins-Strebel direction of $(R,\phi)$} and $(R,\phi)$ is decomposed into cylinders $(V_1,...,V_p)$ where {each $V_j$ has homotopy type $\gamma_j$ and height $b_j$}. 
\end{proposition}
By definition any affine map on $(R,\phi)$ maps all $\phi$-geodesics to $\phi$-geodesics.  
Let $f\in \mathrm{Aff}^+(R,\phi)$ and $D(f)=\bar{A}=\begin{sbmatrix}a&b\\ c&d\end{sbmatrix}\in PSL(2,\mathbb{R})$. Then $f$ maps line segments of direction $\theta\in [0, \pi)$ to line segments of direction $A\theta:=\mathrm{arg}(T_A(e^{\sqrt{-1}\theta}))$. Using the lemma in {\cite[p.56]{A}} we can see that $f$ maps a $\phi$-cylinder of modulus $M$ to a $\phi$-cylinder of modulus $\frac{M}{\sqrt{a^2+c^2}}$. Since the list of moduli of $\phi$-cylinders of one direction are uniquely determined up to order, we have following. 

\begin{lemma}\label{cyl}
Let  $\theta\in J(R,\phi)$ exist and $[M_1^\theta,M_2^\theta,...,M_{n_\theta}^\theta]\in\mathbb{R}_{+}P^{n_\theta-1}$ be the ratio of moduli of the $\phi$-cylinders of direction $\theta$ with $M_j^\theta\geq M_{j+1}^\theta$ $(j=1,2,...,n-1)$.
 If $\bar{A}\in PSL(2,\mathbb{R})$ belongs to $\Gamma (R,\phi)$ then for any direction $\theta\in J(R,\phi)$ following holds. 
\begin{itemize}
\item $A\theta\in J(R,\phi)$
\item $n_{A\theta}=n_\theta$ $(=:n)$
\item $[M_1^{A\theta},M_2^{A\theta},...,M_{n}^{A\theta}]=[M_1^\theta,M_2^\theta,...,M_{n}^\theta]\in\mathbb{R}_{+}P^{n-1}$
\end{itemize}
 \end{lemma}
Next we add an assumption that $\phi$ has two finite Jenkins-Strebel directions $\theta_1, \theta_2\in J(R,\phi)$. We can assume $\theta_1\leq \theta_2$ without loss of generality. In this case, $R$ is obtained by finite collections of parallelograms in the way presented in {\cite[Theorem2]{EG}} (in which we conclude $R$ is finite analytic type even for more general settings). 
We review that construction.

For $i=1,2$ let $\alpha_i=e^{\sqrt{-1}\theta_i}$ and $W_1^i,...,W_{n_i}^i$ be the disjoint $\phi$-cylinders of direction $\theta_i$  which almost every point in $R$ lies on. For each $i,j$ by an analytic continuation of local inverse of $\phi$-charts we construct a holomorphic covering $F^i_j:S^i_j\rightarrow W_j^i$ where $S_j^i=\{0<\mathrm{Im}z<h_j^i\}\subset \mathbb{C}$ and $\mathrm{Deck}(F_j^i)=\langle z\mapsto z+c_j^i \rangle$ for some $h_j^i, c_j^i>0$. 
(Now $c_j^i/h_j^i$ is the modulus of $W_j^i$.) 
By construction $F_j^{i*}(\alpha_i\phi)=d(z_j^i)^2$ holds for any $\phi$-coordinate $z_j^i$ in $W_j^i$. 

\begin{figure}[htbp]
\begin{center}
\includegraphics[width=110mm]{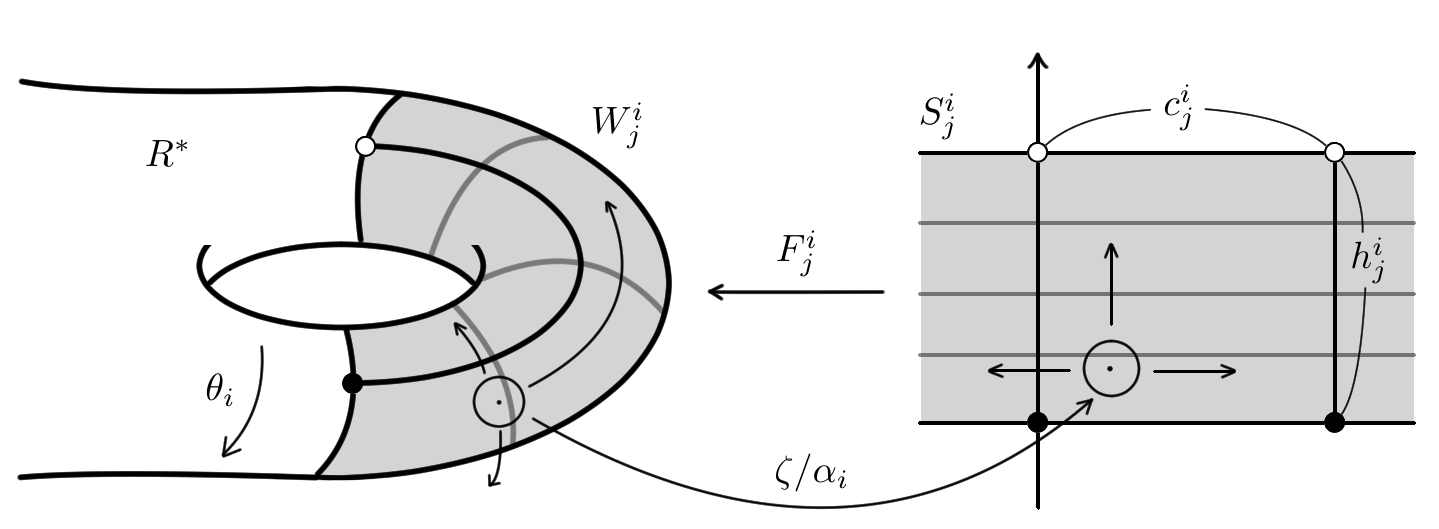}
  \caption{$\phi$-cylinder $W_j^i$ and covering $F^i_j$}
\label{Fig.2}    
\end{center}
\end{figure}

For any $p\in S_j^1$, there is a neighborhood $U$ in which $F_j^1=F_k^2\circ f$ for some $k$ and some holomorphic function $f:U\rightarrow S_k^2$. 
Now $f^*(d(z_j^2)^2)=f^*(F_k^{2*}(\alpha_2\phi))=F_j^{1*}(\alpha_2\phi)=(\alpha_2/\alpha_1)d(z_j^1)^2$ and this implies that $f$ is of the form $f(z_j^1)=\alpha z_j^1+\beta$ where $\alpha\in\{\pm\sqrt{\alpha_2/\alpha_1}= \pm e^{\sqrt{-1}(\theta_2-\theta_1)}\}$ and $\beta\in \mathbb{C}$. 
Thus $U$ is a subset of $V_k=\{z_j^1\in S_j^1|\alpha z_j^1+\beta\in S_k^2\}$. By analytic continuation we see that $F_j^1=F_k^2(\alpha z_j^1+\beta)$ on $V_k$. 
If we replace $p$ by $p+mc_j^1$ for some $m\in\mathbb{Z}$ then $F_j^1=F_k^2(\alpha z_j^1+\beta ')$ still holds on $V_k+mc_j^1$ where $\beta '=\beta-m\alpha c_j^1$. 
So the condition $F_j^1(z_j^1)=F_k^2(\alpha z_j^1+\beta ')$ on $V_k$ for some $\beta '\in\mathbb{C}$ is preserved by covering transformations of $F_j^1$. 

\begin{figure}[htbp]
\begin{center}
\includegraphics[width=110mm]{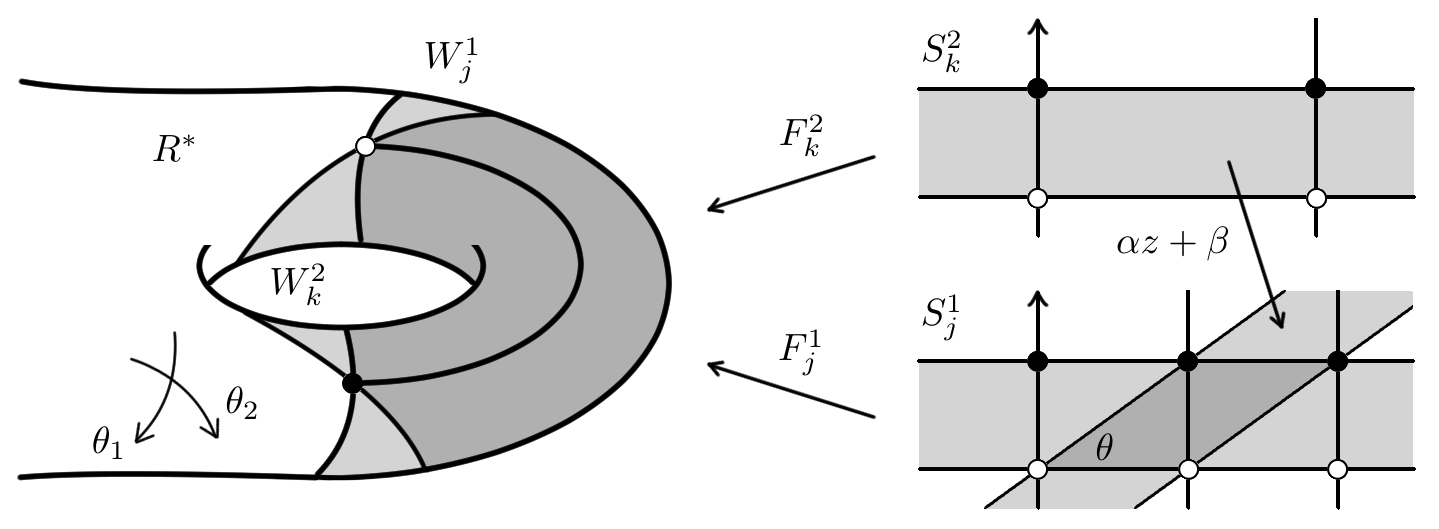}
  \caption{Parallelogram as an intersection of two $\phi$-cylinders}
\label{Fig.3}    
\end{center}
\end{figure}
The parallelogram $V_k$ is mapped into $W_k^2$ and so it does not intersect $V_{k'}$ with $k'\neq k$. Hence the parallelograms on which $F_{j}^1(z_j^1)=F_k^2(\alpha z_j^1 +\beta)$ fill the strip $S_j^1$ by translations in $c_j^1\mathbb{Z}$. Similarly we have same statement for those parallelograms with respect to the strips $S_k^2$. By choosing finite collections of  those parallelograms and gluing them at the points which are mapped to same point by $F_j^1$ or $F_k^2$, we have a compact surface $\tilde{R}$. Finally we obtain $R$ by removing at most finitely many vertices of the parallelograms from $\tilde{R}$, as of finite analytic type. 

We have made a decomposition of $R$ into finite parallelograms $\{V_\lambda\}_{\lambda\in\Lambda}$, where for any $\lambda$ there exists $j_\lambda,k_\lambda$ such that $F^1_{j_{\lambda}}(V_\lambda)=F^2_{k_{\lambda}}(V_\lambda)=W_{j_{\lambda}}^1\cap W_{k_{\lambda}}^2$. 
Now $V_\lambda$ has an angle $\theta=\theta_2-\theta_1$ and a modulus $M_\lambda=(h_{j_{\lambda}}^1/h_{k_{\lambda}}^2) \sin \theta$ for each $\lambda$. 
For each $\theta_1,\theta_2\in J(R,\phi)$ we call these parallelograms the $(\theta_1,\theta_2)$\textit{-parallelograms} on $(R,\phi)$. 
We remark that in above construction each embedding $F_{j_\lambda}^1:V_\lambda\hookrightarrow R$ of $(\theta_1,\theta_2)$-parallelogram is up to inter-decompositions with $z\mapsto \pm z+c$ uniquely determined by $j_\lambda$.

Any $f\in\mathrm{Aff}^+(R,\phi)$ with $D(f)=\bar{A}\in PSL(2,\mathbb{R})$ maps each $\phi$-cylinder of direction $\theta$ to the one of direction $A\theta$. So for each $\theta_1,\theta_2\in J(R,\phi)$ the $(\theta_1,\theta_2)$-parallelograms are mapped to $(A\theta_1,A\theta_2)$-parallelograms. 
On the Euclidian plane, the variation of modulus of $(\theta_1,\theta_2)$-parallelogram under an affine map with derivative $A\in SL(2,\mathbb{Z})$ is described as scalar multiple by $|T_A(e^{\sqrt{-1}\theta_2})|/|T_A(e^{\sqrt{-1}\theta_1})|$. 
The same argument can be said for each of $(\theta_1,\theta_2)$-parallelograms on $(R,\phi)$ and we have following lemma. 

\begin{lemma}\label{paral}
Let $\theta_1,\theta_2\in J(R,\phi)$ and $V$ be a $(\theta_1,\theta_2)$-parallelogram on $(R,\phi)$. Then an affine map $f\in\mathrm{Aff}^+(R,\phi)$ with derivative $A\in PSL(2,\mathbb{R})$ varies the modulus $M(V)$ of $V$ by the multiple of $\rho_{A,\theta_1,\theta_2}=|T_A(e^{\sqrt{-1}\theta_2})|/|T_A(e^{\sqrt{-1}\theta_1})|$. 
In particular, if $(R,\phi)$ is decomposed into $(\theta_1,\theta_2)$-parallelograms $\{V_\lambda\}_{\lambda=1}^N$ then $(M(f(V_1)),M(f(V_2)),...,M(f(V_N)))= \rho_{A,\theta_1,\theta_2}\cdot(M(V_1),M(V_2),...,M(V_N))$ holds. 


\end{lemma}
Furthermore the structure how the edges of those parallelograms are glued is preserved by a homeomorphism $f$. 
In the next section we construct a combinatorial characterization to present this situation more precisely.


\subsection{$F_2-$action to parallelograms}\label{F2}
We continue the assumptions of $(R,\phi)$ and notations as in last section. 
\\
Fix 
signs of embeddings $F_{j_\lambda}^1:V_\lambda\hookrightarrow R$ of $(\theta_1,\theta_2)$-parallelograms for each $j_\lambda$. 
We construct a group $\hat{G}$ which represents how the $(\theta_1,\theta_2)$-parallelograms $\{V_\lambda\}$ are glued. 

For each $\hat{\lambda}=(\lambda,\varepsilon)\in\Lambda\times\{\pm 1\}$, fix some interior point $p_{\hat{\lambda}}\in V_\lambda$ and take a segment $\gamma_{\hat{\lambda}}$ to the direction $\varepsilon$ (resp.\ $\varepsilon e^{\sqrt{-1}\theta}$) from $p_{\hat{\lambda}}$ to the boundary point $p_{\hat{\lambda}} '\in \partial V_\lambda$. We take unique $\lambda '\in \Lambda$ so that $F_{j_\lambda}^1(p_{\hat{\lambda}} ')=F_{j_{\lambda '}}^1(q)$ (resp.\ $F_{k_\lambda}^2(p_{\hat{\lambda}} ')=F_{k_{\lambda '}}^2(q)$) for some $q\in \partial V_{\lambda '}$ and $\varepsilon '\in \{\pm 1\}$ so that $\varepsilon =\varepsilon '$ (resp.\ $F_{k_\lambda}^2\circ ( \left.{F_{k_{\lambda '}}^2}\right|_{U})^{-1}(z)=\varepsilon\varepsilon 'z+c$ on some neighborhood $U$ of $q$ for some $c\in\mathbb{C}$). Further we define $x(\lambda,\varepsilon)$ (resp.\ $y(\lambda,\varepsilon)$) by $(\lambda ',\varepsilon ')$. 
Then we will see that $F_2$, the free group generated by $x,y$ acts on $\hat{\Lambda}:=\Lambda\times\{\pm 1\}$. 


\begin{figure}[htbp]
\begin{center}
\includegraphics[width=120mm]{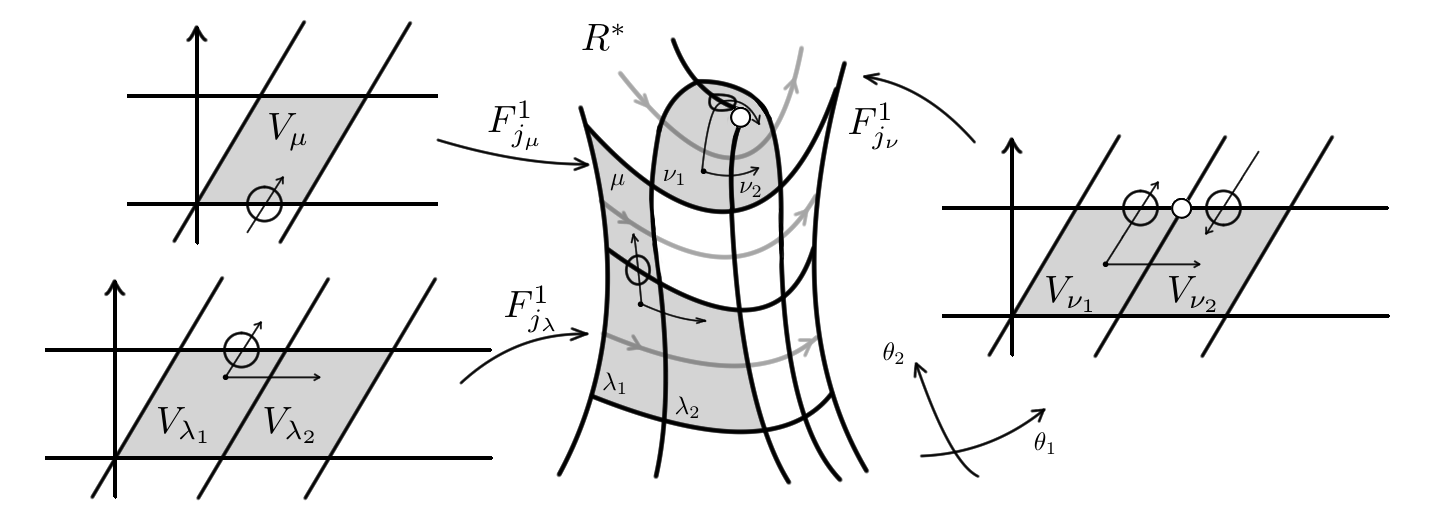}
  \caption{$F_2$-action: In this case we have $\mathbf{x}(\lambda_1,1)=(\lambda_2,1)$, $\mathbf{y}(\lambda_1,1)=(\mu,1)$, $\mathbf{x}(\nu_1,1)=(\nu_2,1)$, and $\mathbf{y}(\lambda_1,1)=(\lambda_2,-1)$. If we choose another sign of embedding of one cylinder of direction $\theta_1$, the signs of elements in that cylinder will be alternated.  }
\label{Fig.4}    
\end{center}
\end{figure}
\begin{lemma}\label{formula}
$x,y:\hat{\Lambda}\rightarrow \hat{\Lambda}$ are bijective. A symmetry $\sigma(\lambda,\varepsilon)=\sigma^{-1}(\lambda,-\varepsilon)$ holds for any $\sigma=x,y$ and $(\lambda,\varepsilon)\in\hat{\Lambda}$. Furthermore, $\phi$ is Abelian if and only if $y$ fixes the signs of elements in $ \hat{\Lambda}$ or equivalently the action $F_2 \hspace{0mm} ^\curvearrowright \hat{\Lambda}$ descends via the projection $p_1:\hat{\Lambda}\rightarrow\Lambda$. 
\begin{pf}
We define $x^{-1},y^{-1}:\hat{\Lambda}\rightarrow \hat{\Lambda}$ as same as $x,y$ except for the direction of $\gamma_{\hat{\lambda}}$ for each ${\hat{\lambda}}$, with it reversed. Then they give inverse maps in the sense of $\lambda\in\Lambda$ at least. For each ${\hat{\lambda}}=(\lambda,\varepsilon)$, $y({\hat{\lambda}})=(\lambda ',\varepsilon_1)$, and $y^{-1}\circ y(\hat{\lambda})=(\lambda,\varepsilon_2)$, $\varepsilon=\varepsilon_1$ if and only if coordinate transformations around the edge between $V_\lambda$ and $V_{\lambda'}$  is of the form $z\mapsto -z+c$, thus it is equivalent to $\varepsilon_2=\varepsilon$. 
So $\varepsilon_1=\varepsilon_2$, $x^{-1},y^{-1}$ give inverse maps indeed. We can see easily from the construction that $\sigma=x,y$ satisfies $\sigma(\lambda,\varepsilon)=\sigma^{-1}(\lambda,-\varepsilon)\cdots(*)$. 

We made $(\theta_1,\theta_2)$-parallelograms $\{V_\lambda\}$ to be glued so that there is no edge in the direction $\theta_1$ where coordinate transformations are of the form $z\mapsto -z+c$. So $\phi$ is Abelian if and only if there is no edge in the direction $\theta_1$ like that, equivalently $y$ fixes the signs of elements in $\hat{\Lambda}$. By the formula $(*)$ this is equivalent to that the action $F_2 \hspace{0mm} ^\curvearrowright \hat{\Lambda}$ descends via $p_1:\hat{\Lambda}\rightarrow\Lambda$. 
\qed\end{pf}

\end{lemma}

We denote the homomorphism which gives the action $F_2 \hspace{0mm} ^\curvearrowright \hat{\Lambda}$ by $\hat{m}:F_2\rightarrow\mathrm{Sym}(\hat{\Lambda})$, $\hat{m}(w)=\hat{m}_w$ for each $w\in F_2$, and $\hat{G}=\langle \hat{m}_x,\hat{m}_y \rangle<\mathrm{Sym}(\hat{\Lambda})$. 
In the case that $\phi$ is Abelian, we also denote ones of the projected action $F_2 \hspace{0mm} ^\curvearrowright \Lambda$ by ${m}:F_2\rightarrow\mathrm{Sym}(\Lambda)$, ${m}(w)={m}_w$, and ${G}=\langle {m}_x,{m}_y \rangle<\mathrm{Sym}(\Lambda)$. 

Let $R^*$ be the surface obtained by puncturing $R$ at all the points which are vertices of the $(\theta_1,\theta_2)$-parallelograms. Let $\hat{R}$ be the translation surface given by taking a double of $R^*$ like non-oriented origamis. $\hat{R}$ is decomposed into the $(\theta_1,\theta_2)$-parallelograms $\{V_{\hat{\lambda}}\}_{\hat{\lambda}\in\hat{\Lambda}}$, with the action of $F_2$ as above giving $\hat{\Lambda}$ as `$\Lambda$' and $\hat{G}$ as `$G$'.  
If $\phi$ is not Abelian $\hat{R}$ is the surface given by an analytic continuation of locally defined Abelian differential $dz$ on $R^*$. 

\begin{proposition}Let $\phi$ be Abelian, $\lambda_0\in\Lambda$, and $p_0\in V_{\lambda_0}$. Then we have following. 
\begin{enumerate}\label{act}
\item $H=H_{{G}}:=\mathrm{Stab}_{F_2}(\lambda_0)<F_2$ is isomorphic to $ \pi_1({R},p_0)$. 
\item There is a 1-1 correspondence between $F_2/H_{{G}}$ and $\Lambda$. 
\item ${G}$ acts transitively on ${\Lambda}$. 
\end{enumerate}
\begin{pf}For any path $C$ in $R^*$, we can take a path $L_C$ which is a product of finite  line segments of direction $\theta_1$ or $\theta_2$ joining neighboring parallelograms so that $C,C'$ are homotopic with fixed endpoints. We define $w_C\in F_2$ to correspond to the order of segments in $L_C$, by replacing segments in the direction $\theta_1,\theta_2$ by $x,y$ respectively. 

For each $w\in F_2$ we take a path $L_{w}$ in $R^*$ starting from $p_0$ to some point in $V_{m_w(\lambda_0)}$, composed of line segments in the direction $\theta_1,\theta_2$ joining neighboring parallelograms whose order respects the one of $x,y$ in $w$. 
$L_w$ is uniquely determined up to variation of the end point in the parallelogram and homotopy with fixed endpoints. 
We define a homomorphism $\Phi:F_2\rightarrow \Lambda$ by $w\mapsto\Phi(w)$ where the end point of $L_w$ belongs to $V_{\Phi(w)}$. Now the kernel of $\Phi$ is $H_G$. 
So we have an isomorphism between $H$ and $\pi_1({R},p_0)$ by taking $L_w$ for each $w\in H$ so that the endpoint of $L_C$ is $p_0$. So (a),(b) follows. 

For any $\lambda,\lambda'\in\Lambda$ there is a path $C$ joining $V_\lambda$ and $V_{\lambda'}$ since $R^*$ is connected. We see $w\in F_2$ has an action defined by the segments of $L_C$, sending $\lambda$ to $\lambda'$. Thus (c) holds. \qed\end{pf}
\end{proposition}



\subsection{Characterization of surfaces}\label{sec3-3}

We consider what condition is needed for parallelograms $\{V_\lambda\}_{\lambda\in\Lambda}$ and a permutation group $\hat{G}=\langle \hat{m}_x,\hat{m}_y \rangle<\mathrm{Sym}(\hat{\Lambda})$ to form a flat surface where they are the total $(\theta_1,\theta_2)$-parallelograms and the one given by the action $F_2 \hspace{0mm} ^\curvearrowright \hat{\Lambda}$ respectively.

First we continue the assumptions for the decomposition of $(R,\phi)$.  
For each parallelogram $V_\circ$, we denote the modulus by $M_\circ$ and the area by $S_\circ$.  
By the construction of $\hat{R}$ the parallelograms $V_\lambda$ in $R^*$ and $V_{(\lambda,1)},V_{(\lambda,-1)}$ in $\hat{R}$ coincide by translation. 
In this sense we define $h_{(\lambda,\cdot)}:=h_\lambda$ for each $\lambda\in\Lambda$. 
We have $M_\lambda=(h_k^2/h_j^1)\sin \theta$, $S_\lambda=(h_j^1h_k^2)\sin \theta$, $h_j^1=h_{\hat{m}_x(j,\varepsilon)}^1$, $h_k^2=h_{\hat{m}_y(k,\varepsilon)}^2$, and in particular following condition for the decomposition of $(R,\phi)$ into $(\theta_1,\theta_2)$-parallelograms. 
\begin{lemma}\label{MS}
For each $\hat{\lambda}\in\hat{\Lambda}$, 
$S_{\hat{m}_{x}(\hat{\lambda})}=\dfrac{M_{\hat{\lambda}}}{M_{\hat{m}_{x}(\hat{\lambda})}}S_{\hat{\lambda}}$ and $S_{\hat{m}_{y}(\hat{\lambda})}=\dfrac{M_{\hat{m}_{y}(\hat{\lambda})}}{M_{\hat{\lambda}}}S_{\hat{\lambda}}$.

\end{lemma}

Conversely we consider $\Lambda=\{1,2,..., N\}$ for $N\in\mathbb{N}$, $\hat{\Lambda}=\Lambda\times\{\pm1\}$, and an arbitrary pair of $M=[M_1,M_2,...,M_N]\in\mathbb{R}_{+}P^{N-1}$ and $\hat{G}=\langle \hat{m}_x,\hat{m}_y \rangle<\mathrm{Sym}(\hat{\Lambda})$. We assume the symmetry of $\hat{G}$ as the formula in Lemma \ref{formula}. 
An $N$-tuple of parallelograms $V=(V_1,V_2,...,V_N)$ with moduli list $M$ is uniquely determined up to congruence by an area list $S=(S_1,S_2...,S_N)\in \mathbb{R}_{+}^N$. 
For each $w\in F_2$ and $\hat{\lambda}\in\hat{\Lambda}$, the formulae in Lemma \ref{MS} defines the unique area $S_{\hat{m}_w(\hat{\lambda})}(\hat{\lambda},w)$ which is necessary for $V_{p_1(\hat{m}_w(\hat{\lambda}))}$ to be glued via the path $L_w$ on the rule given by $\hat{G}$, with a flat structure naturally given. 
For $V$ to form $(R,\phi)$ with the action of $F_2$ given by $\hat{G}$, the area $S_{\hat{m}_w(\hat{\lambda})}(\hat{\lambda},w)$ should depend only on $\lambda:=p_1(\hat{m}_w(\hat{\lambda}))$ and determine $S_{\lambda}$.


Let $\gamma_{-I}$ be the automorphism of $F_2$ defined by $(x,y)\mapsto (x^{-1},y^{-1})$. 
We denote $-(\lambda,\varepsilon):=(\lambda,-\varepsilon)$ for each $(\lambda,\varepsilon)\in\hat{\Lambda}$. Similar to Proposition \ref{origami} the condition for $\hat{G}=\langle \mathbf{x,y}\rangle$ to give an origami which comes from the double cover of some flat surface as before is characterized by conditions for `monodromy map' $\hat{m}$. That is, 
\begin{enumerate}
\item (symmetry) $\hat{m}_{w}(\hat{\lambda})=-\hat{m}_{\gamma_{-I}(w)}(\hat{\lambda})$ for any $\hat{\lambda}\in\hat{\Lambda}$ and $w\in F_2$, 
\item (non-branching) $\mathbf{y}(\hat{\lambda})\neq -\hat{\lambda}$ for any $\hat{\lambda}\in\hat{\Lambda}$, and 
\item (connectivity) the action $\hat{G}\hspace{0mm} ^\curvearrowright \hat{\Lambda}$ is transitive with respect to first ingredients. 
\end{enumerate}

\begin{definition}\label{def}
Let $N\in\mathbb{N}$, $\Lambda=\{1,2,...,N\}$, $\hat{\Lambda}=\Lambda\times\{\pm1\}$, $M=[M_1,M_2,...,$\\$M_N]\in \mathbb{R}_{+}P^{N-1}$, and $\hat{m}:F_2\rightarrow\mathrm{Sym}(\hat{\Lambda})$ be a homomorphism with 
three conditions stated above. 
We define $K_{\mathcal{O}}=K_{M,\hat{G}}:\hat{\Lambda}\times F_2\rightarrow \mathbb{R}$ by following. 
\begin{itemize}
\item $K_{\mathcal{O}}(\ \cdot\ ,1)=1$. 
\item For any $\hat{\lambda}\in\hat{\Lambda}$, $K_{\mathcal{O}}(\hat{\lambda}, x)=\frac{M_{\hat{\lambda}}}{M_{\hat{m}_{x}(\hat{\lambda})}}$ and $K_{\mathcal{O}}(\hat{\lambda}, y)=\frac{M_{\hat{m}_{y}(\hat{\lambda})}}{M_{\hat{\lambda}}}$. 
\item For any $w_1,w_2 \in F_2$ and $\hat{\lambda}\in \hat{\Lambda}$, $K_{\mathcal{O}}(\hat{\lambda},w_1w_2)=K_{\mathcal{O}}(\hat{\lambda},w_1)K_{\mathcal{O}}(\hat{m}_{w_1}(\hat{\lambda}),w_2)$
\end{itemize}
We call $\mathcal{O}=(M,\hat{G}=\langle \mathbf{x,y}\rangle)$ an \textit{extended origami} of degree $N$ if $M\in\mathbb{R}_{+}P^{N-1}$, $\hat{G}<S_{2N}$ with $(\mathbf{x,y})=(\hat{m}(x),\hat{m}(y))$ for some $\hat{m}$ as above, 
 and $K_{\mathcal{O}}(1,w)=1$ for all $w\in H_{\hat{G}}$. Extended origamis $\mathcal{O}_i=(M^i=[M_1^i,M_2^i,...,M_N^i],$\\$G_i=\langle \mathbf{x}_i,\mathbf{y}_i \rangle)$ $(i=1,2)$ of order $N$ are \textit{isomorphic} if there is a pair $(\Phi,\sigma)$ of $\Phi:G_1\rightarrow G_2$ and $\sigma\in S_{2N}$ such that 

\begin{enumerate}
\item $\Phi:\hat{G_1}\rightarrow \hat{G_2}$ is an isomorphism with $(\Phi(\mathbf{x}_1),\Phi(\mathbf{y}_1))=(\mathbf{x}_2,\mathbf{y}_2)$, 
\item $[M_{p_1\circ\sigma(1)}^1, M_{p_1\circ\sigma(2)}^1,..., M_{p_1\circ\sigma(N)}^1]=[M_1^2, M_2^2,..., M_N^2]$, and
\item $\sigma(\hat{m}_w(\hat{\lambda}))=\hat{m}_{\Phi(w)}(\sigma(\hat{\lambda}))$ for each $\hat{\lambda}\in\hat{\Lambda}$, $w\in \hat{G}$. 
\end{enumerate}
We call $(\Phi,\sigma)$ an \textit{isomorphism} between extended origamis $\mathcal{O}_1$ and $\mathcal{O}_2$. 

\end{definition}
By the symmetry of $\hat{G}$, if $\mathbf{g}\in \hat{G}$ contains a cycle $(\hat{\lambda}_1\hat{\lambda}_2...\hat{\lambda}_n)$ then it also contains the cycle $(-\hat{\lambda}_n-\hat{\lambda}_{n-1}...-\hat{\lambda}_1)$. 
From now on we omit half of the cycles in $\mathbf{x},\mathbf{y}\in\hat{G}$ and denote $(\lambda,1),(\lambda,-1)\in\hat{\Lambda}$ by $\lambda , \lambda^-$ respectively.

\begin{theorem}\label{thm1}
A compact flat surface $(R,\phi)$ with a pair of two distinct Jenkins-Strebel directions $(\theta_1,\theta_2)\in J(R,\phi)^2$ is up to isomorphism uniquely determined by a triple $P(R,\phi,(\theta_1,\theta_2))=(\Theta, k,\mathcal{O})$ where $\Theta=(\theta_1,\theta_2)\in[0,\pi)^2$ with $\theta_1\neq\theta_2$, $k>0$, and $\mathcal{O}$ is an extended origami. 
\begin{pf}
Aa we have already seen the decomposition of $(R,\phi)$ into directions $(\theta_1,\theta_2)$ determines 
an extended origami $\mathcal{O}$. We take $k>0$ as the modulus of parallelogram labelled with $1\in\Lambda$. 
 
Conversely if an extended origami $\mathcal{O}=(M=[M_1,M_2,...,M_N],\hat{G}=\langle\mathbf{x,y}\rangle)$ of degree $N\in\mathbb{N}$, $k>0$, and $(\theta_1,\theta_2)\in[0,\pi)^2$ with $\theta_1\neq\theta_2$ are given, then we can construct a flat surface $(R,\phi)$ with $\mathcal{O}(R,\phi,(\theta_1,\theta_2))=\mathcal{O}_1$ as follows. We take an $N$-tuple of Euclidian $(\theta_1,\theta_2)$-parallelograms with the moduli list $k/\Sigma M_\lambda \cdot(M_1,M_2,...,M_N)$. We glue them by the rule given by $\hat{G}$, so that $\mathbf{x,y}$ correspond to each segments of direction $\theta_1,\theta_2$ respectively. Let $R$ be the resulting surface and $R^*$ be surface given by puncturing $R$ at all the vertices. Now natural coordinates $z$ given by those parallelograms (as well origamis) define the quadratic differential $\phi=dz^2$ on $R^*$, for which those parallelograms are $(\theta_1,\theta_2)$-parallelograms. $\phi$ is uniquely extended to $R$. 
It is clear that flat surfaces with isomorphic P-decompositions are isomorphic. 
 
If two extended origamis $\mathcal{O}_1,\mathcal{O}_2$ of same order give isomorphic flat surfaces $(R,\phi),(S,\psi)$ under common $(\theta_1,\theta_2)$ and $k$, then there exists a locally affine quasiconformal homeomorphism $f:R\rightarrow S$ with derivative $[I]$. 
$f$ descends to a map between $(\theta_1,\theta_2)$-parallelograms on $R$ and $S$ of the form $z\mapsto \eta z+c$. This representation can be  extended to each cylinder and $\eta\in\{\pm 1\}$ is unique for each cylinders of direction $\theta_1$. 
We define $\sigma(\lambda,\varepsilon)=(\lambda',\eta\varepsilon)$ where $\lambda$-th parallelogram on $R$ is mapped to $\lambda'$-th parallelogram on $S$, which preserves the ratio of moduli. 
For each $(\lambda,\varepsilon)\in\hat{\Lambda}$ geodesics starting from $\lambda$-th parallelogram used to define $\hat{m}_x(\lambda,\varepsilon)$, $m_y(\lambda,\varepsilon)$ are mapped to ones for $\hat{m}_x(\sigma(\lambda,\varepsilon\eta))$, $\hat{m}_y(\sigma(\lambda,\varepsilon\eta))$ respectively. 
Thus an isomorphism between permutation groups of $\mathcal{O}_1,\mathcal{O}_2$ is induced as $(\mathbf{x_1,y_1})\mapsto (\mathbf{x_2,y_2})$ and to be compatible with $\sigma$, finally we have an isomorphism between $\mathcal{O}_1$ and $\mathcal{O}_2$. 
\qed\end{pf}

\end{theorem}

\begin{definition}\label{Pdec}
We call $P(R,\phi,(\theta_1,\theta_2))=((\theta_1,\theta_2),k,\mathcal{O})$ in Theorem \ref{thm1} a \textit{P-decomposition} of $(R,\phi)$ into directions $(\theta_1,\theta_2)$. 
For a P-decomposition $((\theta_1,\theta_2),k,\mathcal{O})$  and $\bar{A}\in PSL(2,\mathbb{Z})$, we define $\bar{A}\cdot((\theta_1,\theta_2),k,\mathcal{O})$ by the P-decomposition $((A \theta_1,A\theta_2),$\\$\rho_{A,\theta_1,\theta_2}k,\mathcal{O})$. 
We say two P-decompositions $(\Theta_j,k_j,\mathcal{O}_j)$ $(j=1,2)$ are \textit{isomorphic} if $\Theta_1=\Theta_2$, $k_1=k_2$, and $\mathcal{O}_1$ is isomorphic to $\mathcal{O}_2$. 
\end{definition}

\begin{exm}
$\Theta=\Theta_0:=(0,\frac{\pi}{2})$, $k=1$, and an extended origami with $M=\mathbf{1}:=[1,1,..,1]$ give a P-decomposition corresponding to an origami, which is a double covering of original surface.  \label{exm1}
\\Let us consider a flat surface as shown in Fig.\ \ref{Fig.5}. This figure determines a flat surface $(R,\phi)$ of type $(1,6)$.  The orders of $\phi$ are $-1$ at three points, $0$ at two points, and $4$ at one point. 

\begin{figure}[htbp]
\begin{center}
\includegraphics[width=70mm]{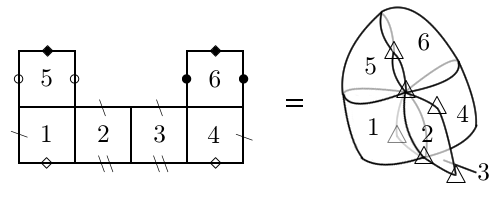}
  \caption{$(R,\phi)$: Edges with the same symbols are glued.}
\label{Fig.5}    
\end{center}
\end{figure}
Now the extended origami comes from $N=6$, $M=\mathbf{1}$, $\mathbf{x}_1=(1234)(5)(6)$, and $\mathbf{y}_1=(156^-4^-)(23^-)$. 
We take the signs of directions of horizontal cylinders in the way shown 
 in Fig.\ref{Fig.6}. 
\begin{figure}[htbp]
\begin{center}
\includegraphics[width=120mm]{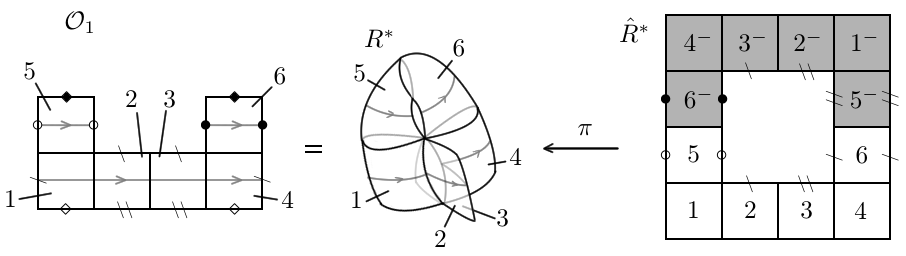}
  \caption{Extended origami $\mathcal{O}_1$: $\mathbf{x}_1=(1234)(5)(6)$, $\mathbf{y}_1=(156^-4^-)(23^-)$}
\label{Fig.6}    
\end{center}
\end{figure}

If we reverse the sign of horizontal cylinder containing the cell labelled with $6$, as in right side of Fig.\ref{Fig.7}, the extended origami has $\mathbf{x}_2=(1234)(5)(6)$ and $\mathbf{y}_2=(1564^-)(23^-)$. It is isomorphic to $(M,\langle\mathbf{x}_1, \mathbf{y}_1\rangle)$ under $6\leftrightarrow 6^-$. 
Similarly if we reverse the signs of all horizontal cylinders, as in left side of Fig.\ref{Fig.7},  we will have an isomorphism given by $\hat{\lambda}\mapsto -\hat{\lambda}$. 

\begin{figure}[htbp]
\begin{center}
\includegraphics[width=125mm]{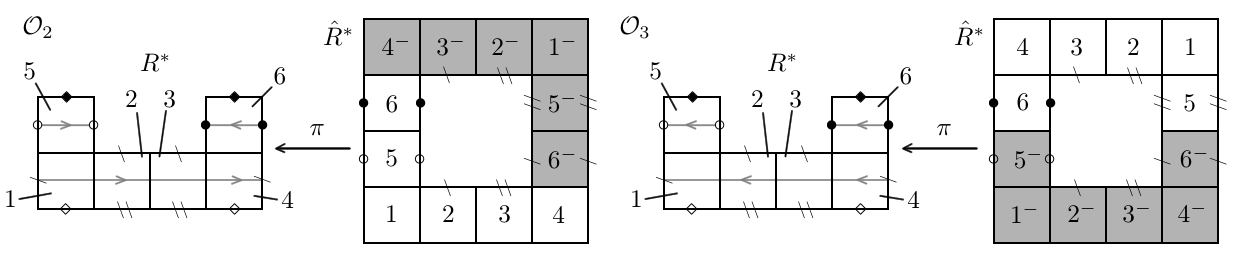}
  \caption{Extended origamis $\mathcal{O}_2$, $\mathcal{O}_3$: Both of them are isomorphic to $\mathcal{O}_1$. }
\label{Fig.7}    
\end{center}
\end{figure}

\end{exm}

An affine map on $(R,\phi)$ with derivative $\bar{A}$ gives a correspondence between 
the initial decomposition $P(R,\phi,(\theta_1,\theta_2))$ and the terminal decomposition $P(R,\phi,(A\theta_1,$\\$A\theta_2))$ as 
we stated in Lemma \ref{cyl} and Lemma \ref{paral}. 
By Theorem \ref{thm1} the existence of such an affine map is described as the possibility of decomposition and the correspondence of those P-decompositions. 
So we conclude as follows. 

\begin{theorem}\label{thm2}
Let $(R,\phi)$ be a flat surface with two distinct Jenkins-Strebel directions $\theta_1,\theta_2\in J(R,\phi)$. $\bar{A}\in PSL(2,\mathbb{R})$ belongs to $\Gamma (R,\phi)$ if and only if $A\theta_1,A\theta_2$ belongs to $ J(R,\phi)$ and $ P(R,\phi,(A\theta_1,A\theta_2))$ is isomorphic to ${A}\cdot P(R,\phi,(\theta_1,\theta_2))$. 
\end{theorem}

\begin{rem}\ 
\begin{enumerate}
\item
Theorem \ref{thm1} also holds for surfaces of finite analytic type with no vertices of $(\theta_1,\theta_2)$-parallelograms contained in $R$. In the proof $R$ is obtained by filling in all the punctures of $R^*$. 
The same can be said for Theorem \ref{thm2}, which implies that if ${A}\in \Gamma (R,\phi)$ then the set of vertices of $(A\theta_1,A\theta_2)$-parallelograms on $(R,\phi)$ should coincide with the one of $(\theta_1,\theta_2)$-parallelograms. 
\item For two P-decompositions $P(R_i,\phi_i,\Theta_i)=(\Theta_i,k_i,\mathcal{O}_i)$ $(i=1,2)$ of flat surfaces, the isomorphism between $\mathcal{O}_1,\mathcal{O}_2$ implies that we can take an affine quasiconformal homeomorphism $g:R_1\rightarrow R_2$ as in the proof of Theorem \ref{thm1}. 
The local derivative of $g$ is the one of affine map on Euclidian plane which maps $\Theta_1$-parallelogram of modulus $k_1$ to $\Theta_2$-parallelogram of modulus $k_2$, which is unique matrix up to signs. $g_*(f):=g\circ f\circ g^{-1}$ defines an isomorphism $g_*:\mathrm{Aff}^+(R_1,\phi_1)\rightarrow \mathrm{Aff}^+(R_2,\phi_2)$ and we will see that $\Gamma(R_2,\phi_2)=A\Gamma(R_1,\phi_1)A^{-1}$ where $A\in PSL(2,\mathbb{R})$ is the derivative of $g$, which is defined as same as ones of affine maps. 
\end{enumerate}
\end{rem}

\begin{exm}
We consider the flat surface in Example \ref{exm1}, which has the P-decomposition $(\Theta_0,1,\mathcal{O}=(\mathbf{1},\langle\mathbf{x,y}\rangle))$ where 
$\mathbf{x}=(1234)(5)(6)$, $\mathbf{y}=(156^-4^-)(23^-)$. 
\\
$\bar{T}=\begin{sbmatrix}1&1\\ 0&1\end{sbmatrix}\in PSL(2,\mathbb{Z})$ gives the extended origami defined by $\mathbf{x}_T=(1234)(5)(6)$, and $\mathbf{y}_T=(12^-3564^-)$, which cannot coincide with $\mathbf{x},\mathbf{y}$ under any permutation of cells. 
On the other hand, we see that $\bar{T}^2=\begin{sbmatrix}1&2\\ 0&1\end{sbmatrix}\in PSL(2,\mathbb{Z})$ gives an extended origami $\mathcal{O}_{T^2}$ isomorphic to $\mathcal{O}$. 
So $\bar{T}\not\in \Gamma (R,\phi)$ and $\bar{T}^2\in \Gamma (R,\phi)$. 
\end{exm}
\begin{figure}[htbp]
\begin{center}
\includegraphics[width=120mm]{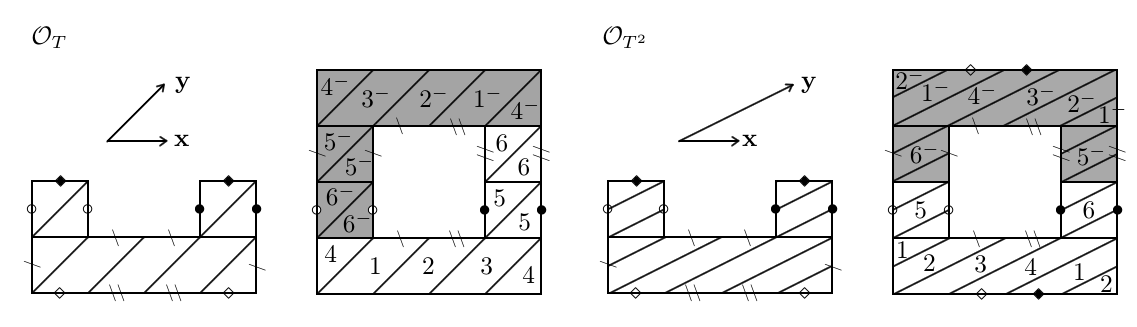}
  \caption{$\mathcal{O}_T$ is not isomorphic to $\mathcal{O}$ and $\mathcal{O}_{T^2}$ is isomorphic to $\mathcal{O}$}
\label{Fig.8}    
\end{center}
\end{figure}

\begin{exm}
We show examples whose moduli lists are not rational. 
\\
\begin{enumerate}\item
For $n\geq 2$, a translation surface coming from a regular $2n$-gon $P_{2n}$ (and its translation coverings) is studied in \cite{EG} and \cite{Sh}. It is obtained by gluing $P_{2n}$ at all opposite edges. Its Veech groups is known to be the group generated by $T_{2n}=\begin{spmatrix}1&\cot{\frac{\pi}{2n}} \\ 0&1\end{spmatrix}$ and $R_{2n}=\begin{spmatrix}\cos{\frac{\pi}{n}} &-\sin{\frac{\pi}{n}} \\ \sin{\frac{\pi}{n}}&\cos{\frac{\pi}{n}}\end{spmatrix}$. 
\\This can be seen as in Fig.\ref{Fig.9} for $n=4$. Now the moduli ratio comes from $[1,1,1,\sqrt{2}, \sqrt{2},\frac{1}{\sqrt{2}}, \frac{1}{\sqrt{2}}]$. 
\begin{figure}[htbp]
\begin{center}
\includegraphics[width=110mm]{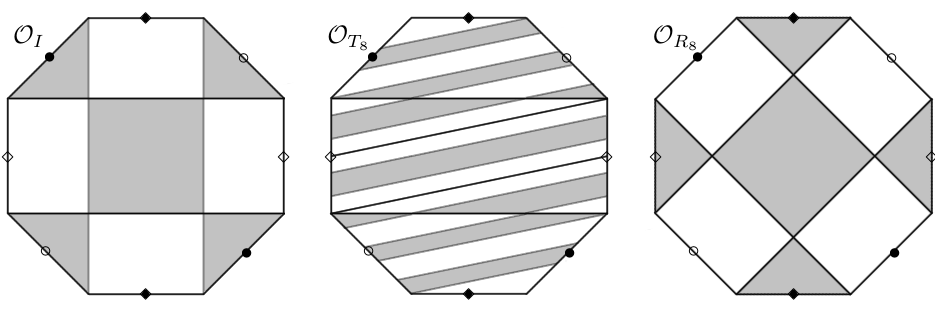}
  \caption{Translation surface coming from $P_8$ and its P-decompositions}
\label{Fig.9}    
\end{center}
\end{figure}

\item
We next consider a flat surface obtained from $P_{2n}$ in the following way. We divide all edges of $P_{2n}$ in half and glue each of them to the adjacent edge in the same direction. As in Fig.\ref{Fig.10} we can do the same observation for decompositions as (a) and so the Veech group includes $\langle T_{2n},R_{2n}\rangle$. 
\begin{figure}[htbp]
\begin{center}
\includegraphics[width=110mm]{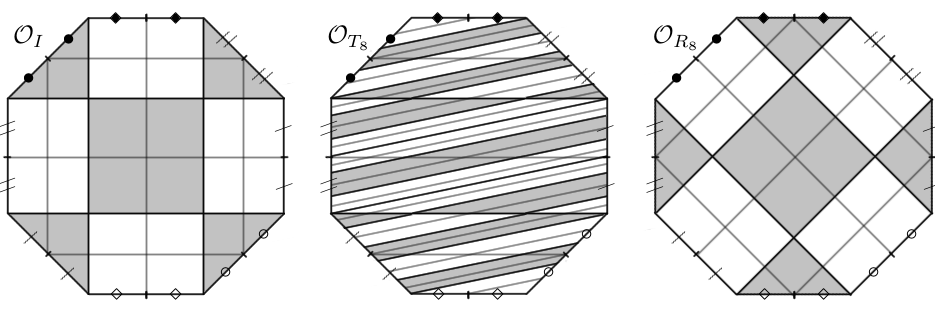}
  \caption{Flat surface coming from $P_8$ and its P-decompositions}
\label{Fig.10}    
\end{center}
\end{figure}

Furthermore, any affine map lifts via the canonical double covering, which is a translation covering of surface in (a). With results in \cite{Sh} we see that the Veech group is contained in $\langle T_{2n},R_{2n}\rangle$, and hence equals to it. 
\end{enumerate}
\end{exm}

\subsection{Further observations}\label{sec3-4}

Next we consider a finite set $V$ of additional marked points in $R$. We assume $R$ is already punctured at all points in $V$. We take the decomposition of $R$ into $(\theta_1,\theta_2)$-parallelograms and obtain $\mathcal{O}=(M,\hat{G})$. 

As same as origamis (see e.g.\ \cite{H}), there is a 1-1 correspondence between cycles in $\hat{m}_{xyx^{-1}y^{-1}}\in \hat{G}$ and vertices of $(\theta_1,\theta_2)$-parallelograms in $\hat{R}$. 
We call a cycle in $\hat{m}_{xyx^{-1}y^{-1}}\in \hat{G}$ an \textit{$\hat{\mathcal{O}}$-vertex}. Two $\hat{\mathcal{O}}$-vertices $v_1,v_2$ correspond to the same point in $R$ if and only if they coincide under the double cover $\hat{R}\rightarrow R^*$, that is $v_2$ equal to $v_1$ or its sign inversion. 
We say such $\hat{\mathcal{O}}$-vertices $v_1,v_2$ are \textit{conjugate}, and call the conjugacy class of an $\hat{\mathcal{O}}$-vertex an \textit{$\mathcal{O}$-vertex}. 
For each $\sigma\in S_{2N}$ and $\hat{\mathcal{O}}$-vertex $v$, we denote by $\sigma^*v$ the cycle obtained by applying $\sigma$ to each element in $v$. 
We call a pair $(\mathcal{O},\mathcal{V})$ of an extended origami $\mathcal{O}$ and a set $\mathcal{V}$ of $\mathcal{O}$-vertices a \textit{marked} extended origami. 

For $(R,\phi)$, $\theta_1,\theta_2\in J(R,\phi)$, and a set $V$ of finite marked points in $R$ we have a P-decomposition $P(R,\phi,(\theta_1,\theta_2),V):=(\Theta, k, (\mathcal{O},\mathcal{V}))$ where $(\Theta,k,\mathcal{O})=P(R,\phi,(\theta_1,\theta_2))$ and $\mathcal{V}$ is the set of $\mathcal{O}$-vertices corresponding to $V$. 

\begin{rem}
Theorem \ref{thm1} can be extended to the general cases of flat surface of finite analytic type by replacing extended origamis with marked extended origamis. In this sense we also call $P(R,\phi,(\theta_1,\theta_2),V)$ a P-decomposition (with marked points). We define isomorphism of marked extended origamis next. 
\end{rem}

Consider an affine map $f\in\mathrm{Aff}^+(R,\phi)$ with derivative $\bar{A}$.  
We have $A\theta_1,A\theta_2\in J(R,\phi)$ and let $(N,\varphi,(\mathcal{O}_A,\mathcal{V}_A)):=P(R,\phi,(A\theta_1,A\theta_2),V)$. 
By Theorem \ref{thm2} there is an isomorphism $(\Phi,\sigma)$ between $\mathcal{O},\mathcal{O}_A$ By the construction given in the proof, $\sigma$ represents how $f$ maps each $(\theta_1,\theta_2)$-parallelograms to $(A\theta_1,A\theta_2)$-parallelograms. For each $[v]\in \mathcal{V}$ by definition \ref{def} (c) $[\sigma^*v]\in\mathcal{V}_A$ and it corresponds to the image of the point in $V$ corresponding to $[v]$. 

We say two marked extended origamis $(\mathcal{O}_i,\mathcal{V}_i)$ $(i=1,2)$ are \textit{isomorphic} if there is an isomorphism $(\Phi,\sigma)$ between $\mathcal{O}_1,\mathcal{O}_2$ and $[v]\mapsto [\sigma^*v]$ gives a well-defined bijection $\mathcal{V}_1\rightarrow \mathcal{V}_2$. 
In such a case for each $[v]\in \mathcal{V}$ the number of elements in $[v],[\sigma^*v]$ coincide and the corresponding vertices have the same valency even in $R$. So permutations among marked points and critical points in $(R,\phi)$ appear at most in the classes of same valencies.

With fixed $\theta_1,\theta_2\in J(R,\phi)$, an affine map on a flat surface $(R,\phi)$ of finite analytic type is characterized as one of $(R^*,\phi)$ extended to stabilize $\partial R\subset \partial R^*$ setwise. 
So we have following. 
\begin{theorem}\label{thm3}
Let $(R,\phi)$ be a flat surface of finite analytic type with two distinct Jenkins-Strebel directions $\theta_1,\theta_2\in J(R,\phi)$. $\bar{A}\in PSL(2,\mathbb{R})$ belongs to $\Gamma (R,\phi)$ if and only if $A\theta_1,A\theta_2$ belongs to $ J(R,\phi)$ and $A\cdot P(R,\phi,(\theta_1,\theta_2),\partial R)$ is isomorphic to $P(R,\phi,(A\theta_1,A\theta_2),\partial R)$. 
\end{theorem}
\begin{corollary}\label{cor}
Let $R$ be a Riemann surface of finite analytic type  and $\phi\in A(\bar{R})$. If for any $p\in \partial R$ all critical points of $\phi$ of order $\mathrm{ord}_p(\phi)$ are contained in $\partial R$, then $\Gamma(R,\phi)=\Gamma(\bar{R},\phi)$. In particular, if there exists a P-decomposition of $(R,\phi)$ whose moduli ratio is rational then $(R,\phi)$ induces a Teichm\"uller curve which is a Belyi surface. 
\begin{pf}
The former claim follows from above observations immediately. 

By taking a set $W$ of sufficiently many additional marked points one can obtain a P-decomposition with all parallelograms congruent. 
Up to conjugation in affine deformations we may assume that $P(R,\phi,\Theta_0,\partial R\cup W)=(\Theta_0,1,(\mathbf{1},\cdot))$. Since developed images in the plane of $\partial R, W$ are infinite sets contained in $\mathbb{Z}+\sqrt{-1}\mathbb{Z}$, we have $\Gamma(R\setminus W,\phi)<\Gamma(R,\phi)<PSL(2,\mathbb{Z})$. (see \cite[Propositioin 2.6]{S1} for details.)

For any $A\in PSL(2,\mathbb{Z})$, $P(R,\phi,A\Theta_0,\partial R\cup W)$ is of the form $(\Theta_0,1,(\mathbf{1},\cdot))$ again. Now the number of cells equals to the one of  $P(R,\phi,\Theta_0,\partial R\cup W)$ and the number of such decompositions are finite. 
So $\Gamma(R\setminus W,\phi)$ is a subgroup of $PSL(2,\mathbb{Z})$ of finite index and so is $\Gamma(R,\phi)$. 
We have conclusion. 
\qed\end{pf}
\end{corollary}

\begin{exm}\label{exmver}
Let us consider the flat surface with $4$ distinguished marked points $a,b,c,d$ as shown in the left of Fig.\ref{Fig.11}.  
\begin{figure}[htbp]
\begin{center}
\includegraphics[width=90mm]{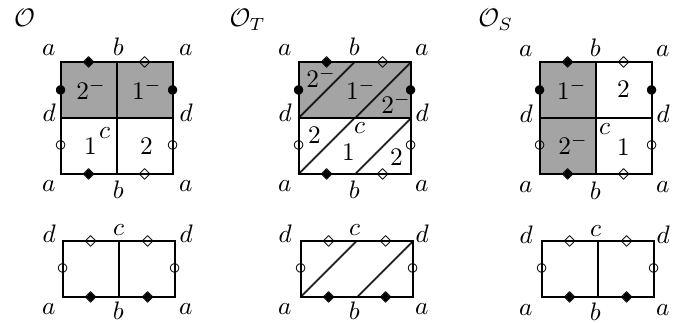}
  \caption{P-decompositions with distinguished marked points}
\label{Fig.11}    
\end{center}
\end{figure}
If we decomposed it into $\Theta_0$, the extended origami $\mathcal{O}=(\mathbf{1},\langle\mathbf{x,y}\rangle)$ is given by $\mathbf{x}=(12)$ and $\mathbf{y}=(12^-)$. 
Now $\mathbf{xyx}^{-1}\mathbf{y}^{-1} =(1)(2)(1^-)(2^-)$ and points $a,b,c,d$ correspond to cycles $(1^-),(2^-),(1),(2)$ respectively. (For instance a path $xyx^{-1}y^{-1}$ from cell $1$ goes around the vertex $c$.)

For $\bar{T}=\begin{sbmatrix}1&1\\ 1&0\end{sbmatrix}$, $\bar{S}=\begin{sbmatrix}0&-1\\ 1&0\end{sbmatrix}\in PSL(2,\mathbb{Z})$, decompositions into $T\Theta_0$,  $S\Theta_0$ gives extended origamis $\mathcal{O}_T$, $\mathcal{O}_S$ which are isomorphic to $\mathcal{O}$. On the other hand the correspondences between vertices and cycles are described as $(a,b,c,d)\leftrightarrow((1^-),(2^-),(2),(1))$ in $\mathcal{O}_T$ and $(a,b,c,d)\leftrightarrow((1^-),(2),(1),(2^-))$ in $\mathcal{O}_S$. They cannot be mapped by any isomorphism of extended origami each other and thus they differ as marked extended origamis. 
Further we can see that the situations in $\mathcal{O}_{T^2}$, $\mathcal{O}_{TS}$, $\mathcal{O}_{ST}$ coincide with $\mathcal{O}$, $\mathcal{O}_{T}$, $\mathcal{O}_{S}$ respectively. 
\end{exm}

\begin{acknowledgements}
I would like to thank Prof.\ Toshiyuki Sugawa for his helpful advices and comments. 
I am grateful to Prof.\ Hiroshige Shiga for his thoughtful guidance with my master thesis, which is a predecessor of this paper. I thank to Prof. Rintaro Ohno for several suggestions. 
Some proposals given by Prof.\ Yoshihiko Shinomiya helped me to get an idea for this paper. 
\\
In Weihnachtsworkshop 2019 at Karlsruhe, I had a lot of significant discussions on my research. I would like to thank Prof.\ Frank Herrlich, Prof.\ Gabriela Weitze-Schmith\"usen, and Prof.\ Martin M\"oller for their expert advices and comments. 
I thank Sven Caspart for helpful discussions. 
\\
\end{acknowledgements}

%
%

\end{document}